\documentstyle[12pt]{article}
\def\R{\ifmmode{\rm I\mkern-3.1mu
R\mkern1mu}\else{\rm I\kern-.18em  R\hskip1pt\ 
}\fi\relax}  
\font\bigreek=Symbol at 16pt
\font\smgreek=Symbol at 10pt
 
\def\b{\beta}

\def\l{\lambda}

\def\n{\nu}

\def\s{\sigma}

\def\so{\underline} 

\def\f{\rightarrow}
\def\q{\forall}
\def\e{\exists}

\def\p{\succ}

\def\R{\ifmmode{\rm I\mkern-3.1mu
R\mkern1mu}\else{\rm I\kern-.18em 
R\hskip1pt\ }\fi\relax} 

\def\Z{\ifmmode{ Z\mkern-4.6mu
Z\mkern2mu}\else{ Z\kern-.28em 
Z\hskip1pt\ }\fi\relax} 

\def\Q{\ifmmode{\rm Q\mkern-10mu
l\mkern4.5mu}\else{\rm Q\kern-.57em
l\hskip3pt\ }\fi\relax} 

\def\N{\ifmmode{\rm I\mkern-3.1mu
N\mkern0.5mu}\else{\rm I\kern-.16em
N\hskip0.5pt\ }\fi\relax} 

\def\C{\ifmmode{\rm C\mkern-8.8mu
l\mkern4mu}\else{\rm C\kern-.48em
l\hskip2.6pt\ }\fi\relax} 

\def\mats{\ifmmode{ {\hbox{\bigreek s}} }\else{ 
{\bigreek s} }\fi\relax}
\def\matsin{\ifmmode{ {\hbox{\smgreek s}} }\else{ 
{\smgreek s} }\fi\relax}
\def\matt{\ifmmode{ {\hbox{\bigreek t}} }\else{ 
{\bigreek t} }\fi\relax}
\def\mattin{\ifmmode{ {\hbox{\smgreek t}} }\else{ 
{\smgreek t} }\fi\relax}

\begin{document}

\vspace*{0cm} 

\begin{center}

{\large \bf An example of a non adequate numeral system}\\ [0,5cm]

{\bf  Karim NOUR} \\ 

\end{center}

{\small {\bf Abstract} {\it A numeral system is defined by three closed
$\l$-terms : a normal $\l$-term $d_0$ for Zero, a $\l$-term $S_d$ for Successor, and
a $\l$-term for Zero Test, such that the $\l$-terms $({S_d}^{i} ~ d_0)$ are
normalizable and have different normal forms. A numeral system is said 
adequate iff it has a closed $\l$-term for Predecessor. This Note gives a simple example of a
non adequate numeral system.}} 

\begin{center}
{\large \bf Un exemple d'un syst\`eme num\'erique non ad\'equat}\\ [0,5cm]

\end{center}

{\small {\bf R\'esum\'e} {\it Un syst\`eme num\'erique est d\'efini par la donn\'ee de trois
$\l$-termes clos: un $\l$-terme normal $d_0$ pour Z\'ero, un $\l$-terme $S_d$ pour le
Successeur, et un $\l$-terme pour le Test \`a Z\'ero, tels que les $\l$-termes $({S_d}^{i} ~
d_0)$ sont normalisables et poss\`edent des formes normales differentes. Un syst\`eme
num\'erique est dit ad\'equat ssi il poss\`ede un $\l$-terme clos  pour le Pr\'ed\'ecesseur.
Dans cette Note nous pr\'esentons un exemple simple d'un syst\`eme num\'erique non
ad\'equat.}}\\

{\large \bf Version Fran\c{c}aise Abr\'eg\'ee} \\

Un {\it syst\`eme num\'erique} est une suite $\bf d$ \rm = $d_0 , d_1 ,..., d_n , ...$
de $\l$-termes normaux clos diff\'erents pour laquelle il existe des $\l$-termes clos $S_d$
et $Z_d$ tels que : \begin{center}
$(S_d ~ d_n) \simeq\sb{\b} d_{n+1}$ pour tout $n \in \N$\\
et\\
$(Z_d ~ d_0) \simeq\sb{\b} \l x \l y x$ \\
$(Z_d  ~ d_{n+1}) \simeq\sb{\b}\l x \l y y$ pour tout $n \in \N$
\end{center}

Les $\l$-termes $S_d$ et $Z_d$ sont appel\'es {\it Successeur} et {\it Test \`a Z\'ero} pour 
$\bf d$. \\
 
Chaque syst\`eme num\'erique peut \^etre naturellement consid\'erer comme un codage des
entiers en $\l$-calcul et donc nous pouvons repr\'esenter les fonctions num\'eriques totales
de la mani\`ere suivante: \\

Une fonction num\'erique totale $\phi : \N^p \f \N$ est dite {\it $\l$-d\'efinissable}
dans un syst\`eme num\'erique {\bf d} si et seulement si :
\begin{center} 
$\e$ $F_{\phi}$ $\q$
$n_1,...,n_p \in \N$ $(F_{\phi} ~ d_{n_1}... d_{n_p}) \simeq\sb{\b} d_{\phi(n_1,...,n_p)}$ 
\end{center}

La diff\'erence entre la d\'efinition d'un syst\`eme num\'erique que nous proposons ici
et celle donn\'ee par H. Barendregt dans [1] est le fait que nous exigons que les $d_i$
soient normaux et diff\'erents. En effet ces derni\`eres conditions permettent avec des
strat\'egies de r\'eduction fix\'ees une fois pour toute (par exemple la strat\'egie de la
r\'eduction gauche qui consiste \`a r\'eduire toujours dans un $\l$-terme le redex le plus \`a
gauche) la valeur exacte d'une fonction calcul\'ee sur des arguments.\\

Un syst\`eme num\'erique {\bf d} est dit {\it ad\'equat} si et seulement s'il existe un
$\l$-terme clos $P_d$ tel que :
\begin{center} 
$(P_d ~ d_{n+1}) \simeq\sb{\b} d_n$ pour tout $n \in \N$ 
\end{center}

Le $\l$-terme $P_d$ est appel\'e {\it Pr\'ed\'ecesseur} pour {\bf d}.  \\

H. Barendregt a d\'emontr\'e dans [1] que :\\

{\it Un syst\`eme num\'erique {\bf d} est ad\'equat si et seulement si toutes les
fonctions r\'ecursives totales sont $\l$-d\'efinissables dans {\bf d}}. \\  

Une question se pose : {\bf Peut on trouver un syst\`eme num\'erique non ad\'equat ?} \\

Nous pr\'esentons dans cette Note un exemple d'un syst\`eme num\'erique non ad\'equat. \\

Le syst\`eme num\'erique non ad\'equat que nous proposons est le suivant :

\begin{center} 
$d_0 = \l x (x ~~ \l x \l y x ~~ \l x x)$ \\
et \\
$d_{n+1} = \l x (x ~~ \l x \l y y ~~ \l x_1 ... \l x_n \l x x)$ pour tout $n \in \N$ 
\end{center}

La d\'emonstration de non ad\'equation s'inspire des techniques d\'evelopp\'ees par J-L
Krivine dans [3] pour d\'emontrer son th\'eor\`eme de mise en m\'emoire.\\

\so{\quad \quad \quad \quad \quad \quad \quad \quad \quad}

\section{Notations and definitions}

The notations are standard (see [1] and [2]).

\begin{itemize}
\item The $\b$-equivalence relation is denoted by $M \simeq\sb{\b} N$. 
\item We denote by $T$ (for True) the $\l$-term $\l x \l y x$ and by $F$ (for False) the
$\l$-term $\l x \l y y$. \item The notation $\s(M)$ represents the result of the simultaneous
substitution $\s$ to the free variables of $M$ after a suitable renaming of the bounded
variables of $M$.     
\item The pair $<M,N>$ denotes the $\l$-term $\l x (x ~ M ~ N)$. 
\item Let us recall that a $\l$-term $M$ either has a {\it head redex} [i.e. $M=\l x_1 ...\l
x_n (\l x U ~ V ~ V_1 ... V_m)$, the head redex being $(\l x U ~ V)$], or is in {\it head
normal form} [i.e. $M=\l x_1 ...\l x_n (x ~ V_1 ... V_m)$].
\item The notation $M \p N$ means that $N$ is obtained from $M$ by some head reductions and we
denote by $h(M,N)$ the length of the head reduction between $M$ and $N$. 
\item A $\l$-term is said {\it solvable} iff its head reduction terminates.
\end{itemize} 

The following results are well known  :

\begin{itemize}
\item[] {\it - If $M$ is $\b$-equivalent to a head normal form then $M$ is solvable.}
\item[]Ê{\it - If $M \p N$, then, for any substitution $\s$, $\s(M) \p \s(N)$,
and $h(\s(M),\s(N))$=h(M,N). In particular, if for some  substitution $\s$, $\s(M)$ is
solvable, then $M$ is solvable.}
\end{itemize}  
 
\section{Numeral systems}

\begin{itemize}
\item A {\it numeral system} is a sequence $\bf d$ \rm = $d_0 , d_1 ,..., d_n , ...$
consisting of different closed normal $\l$-terms such that for some closed $\l$-terms $S_d$
and $Z_d$ : \begin{center}
$(S_d ~ d_n) \simeq\sb{\b} d_{n+1}$ for all $n \in \N$ \\
and\\
$(Z_d ~ d_0) \simeq\sb{\b} T$ \\
$(Z_d  ~ d_{n+1}) \simeq\sb{\b} F$ for all $n \in \N$ 
\end{center}

The $\l$-terms $S_d$ and $Z_d$ are called {\it Successor} and {\it Zero Test} for  $\bf
d$. 
\end{itemize} 
Each numeral system can be naturally considered as a coding of integers in $\l$-calculus
and then we can represent total numeric functions as follows:
\begin{itemize}
\item A total numeric function $\phi : \N^p \f \N$ is {\it $\l$-definable} with respect to a 
numeral  system {\bf d} iff \begin{center}
$\e$ $F_{\phi}$ $\q$ $n_1,...,n_p \in \N$ $(F_{\phi} ~ d_{n_1}... d_{n_p}) \simeq\sb{\b}
d_{\phi(n_1,...,n_p)}$ 
\end{center}
\item A numeral system {\bf d} is called {\it adequate} iff there is a closed
$\l$-term $P_d$ such that
\begin{center} 
$(P_d ~ d_{n+1}) \simeq\sb{\b} d_n$ for all $n
\in \N$. 
\end{center}
The $\l$-term $P_d$ is called {\it Predecessor} for {\bf d}.  
\end{itemize}
H. Barendregt has shown in [1] that :
\begin{itemize}
\item[] {\it A numeral system {\bf d} \it is adequate iff all total recursive functions are
$\l$-definable with respect to} {\bf d}. \\  
\end{itemize}

A question arises : {\bf Can we find a non adequate numeral system ?}

\section{An example of a non adequate numeral systems}

{\bf Theorem} {\it There is a non adequate numeral system.}  \\

\bf Proof \rm For every $n \in \N$, let $p_n = \l x_1...\l x_n \l xx$.\\
Let $d_0 = < T , p_0 >$ and for every $n \geq 1$,  $d_n = <  F , p_n >$.\\ 
It is easy to check that the $\l$-terms $S_d = \l n < F , \l x n >$ and $Z_d = \l n (n ~
T)$ are $\l$-terms for Successor and Zero Test for {\bf d}.\\
Let $\n,x,y$ be different variables. \\
If {\bf d} possesses a closed $\l$-term $P_d$ for Predecessor, then the $\l$-term \\
$Q_d = \l n ((P_d ~ < F , n >) ~ T)$ satisfies the following : 
 \begin{center}
$(Q_d ~ p_{n+1} ~ x ~ y) \simeq\sb{\b} \cases { x &if $n=0$ \cr y &if $n \geq 1$
\cr}$ 
\end{center}
then
 \begin{center}
$(Q_d ~ p_{n+1} ~ x ~ y) \p \cases { x &if $n=0$ \cr y &if $n \geq 1$
\cr}$. 
\end{center}
Therefore  $(Q_d ~ \n ~ x ~ y)$ is solvable and its head normal form does not begin by
$\l$. \\
We have three cases to see :  
\begin{itemize} 
\item   $(Q_d ~ \n ~ x ~ y) \p (x ~ u_1...u_k)$, then $(Q_d ~ p_2 ~ x ~ y) \not \p y$.
\item   $(Q_d ~ \n ~ x ~ y) \p (y ~ u_1...u_k)$, then $(Q_d ~ p_1 ~ x ~ y) \not \p x$. 
\item   $(Q_d ~ \n ~ x ~ y) \p (\n ~ u_1...u_k)$, then $(Q_d ~ p_{k+2} ~ x ~ y) \not \p y$.
\end{itemize}
Each case is impossible. Therefore $\bf d$ \rm is a non adequate numeral system. $\Box$ \\

\bf Acknowledgement.\/ \rm We wish to thank Mariangiola Dezani for helpful discussions.

LAMA - \'EQUIPE DE LOGIQUE \\ 
 \quad \quad UNIVERSIT\'E DE CHAMB\'ERY \\
  \quad \quad \quad \quad 73376 LE BOURGET DU LAC \\

{\it E-mail}: nour@univ-savoie.fr

\end{document}